\documentclass[12pt,centertags,oneside]{amsart}
\usepackage[a4paper,margin=2.5cm]{geometry}
\usepackage[pagebackref]{hyperref}
\usepackage{amsmath,amssymb,amsfonts,amsthm}
\usepackage[abbrev]{amsrefs}
\usepackage[shortlabels]{enumitem}
\usepackage{mathrsfs,amsbsy,bm,mathtools,color}

\newtheorem{theorem}{Theorem}[section]

\newtheorem{prop}[theorem]{Proposition}
\newtheorem{corollary}[theorem]{Corollary}
\newtheorem{lemma}[theorem]{Lemma}
\theoremstyle{definition}
\newtheorem{defn}[theorem]{Definition}
\newtheorem{definition}[theorem]{Definition}
\newtheorem{remark}[theorem]{Remark}

\newtheorem{example}[theorem]{Example}

\newtheorem{question}[theorem]{Question}

\newcommand{\bC}{\mathbb{C}}

\newcommand{\bQ}{\mathbb{Q}}

\newcommand{\mleb}{\mathfrak{m}}

\newcommand{\Q}{\mathbb{Q}}
\newcommand{\CC}{\mathbb{C}}
\newcommand{\R}{\mathbb{R}}
\newcommand{\Z}{\mathbb{Z}}

\newcommand{\lden}{\underline{\sf ud}}
\newcommand{\uden}{\overline{\sf ud}}

\renewcommand{\le}{\leqslant}
\renewcommand{\ge}{\geqslant}

\begin{document}

\title{On a question of Astorg and Boc Thaler}

\author{Zhangchi Chen} 
\address[Zhangchi Chen]{School of Mathematical Sciences, Key Laboratory of MEA (Ministry of Education) and Shanghai Key Laboratory of PMMP, East China Normal University, Shanghai 200241, China}
\email{zcchen@math.ecnu.edu.cn}

\author{Zihao Ye}
\address[Zihao Ye]{School of Mathematical Sciences, East China Normal University, Shanghai 200241, China}
\email{51265500041@stu.ecnu.edu.cn}

\author{Weizhe Zheng}
\address[Weizhe Zheng]{Morningside Center of Mathematics, Academy of Mathematics and Systems 
Science, Chinese Academy of Sciences, Beijing 100190, China; University of 
Chinese Academy of Sciences, Beijing 100049, China} \email{wzheng@math.ac.cn} 

\keywords{Wandering domain, Pisot number, Linear recurrence, Weyl's 
equidistribution theorem}

\subjclass{Primary 37F10; Secondary 11J71, 11K16, 37F80}

\begin{abstract}
Astorg and Boc Thaler studied the dynamics of certain skew-products $f$ 
tangent to the identity on $\CC^2$, with two real parameters $\alpha>1$ and 
$\beta$ derived from its coefficients. They proved that if there exists a 
strictly increasing sequence of positive integers $(n_k)_{k\geqslant 1}$ 
such that $(\sigma_k)_{k\geqslant 1}:=(n_{k+1}-\alpha n_k-\beta\ln 
n_k)_{k\geqslant 1}$ converges, then  $f$ admits wandering domains of rank 
one. They also proved that for $\alpha>1$ with the Pisot property, the 
condition that $\theta:=\frac{\beta\ln\alpha}{\alpha-1}$ is rational is 
sufficient for the existence of $(n_k)_{k\geqslant 1}$ such that 
$(\sigma_k)_{k\geqslant 1}$ converges to a cycle. They asked if this 
condition is necessary. 

When $\alpha$ is an algebraic number, we answer the question of Astorg and 
Boc Thaler in the affirmative. Furthermore, denoting by 
$P(x)\in\mathbb{Z}[x]$ the minimal polynomial of~$\alpha$, we prove that 
$\theta\in\frac{1}{P(1)}\mathbb{Z}$ is necessary and sufficient for the 
existence of $(n_k)_{k\geqslant 1}$ such that $(\sigma_k)_{k\geqslant 1}$ 
converges. Combined with the work of Astorg and Boc Thaler, our result 
provides explicit new examples of skew-products on $\mathbb{C}^2$ with 
wandering domains of rank one.\end{abstract} \maketitle

\section{Introduction}
The existence of wandering domains is a fascinating problem in complex 
dynamics. No rational self-maps of $\mathbb{CP}^1$ have wandering 
domains~\cite{Sullivan-1985}. However, certain algebraic self-maps of 
$\bC^2$ have wandering domains of rank zero~\cite{ABDPR-2016}. 
In~\cite{Astorg-Boc Thaler-2024}, Astorg and Boc Thaler studied the 
dynamics of holomorphic self-maps of $\bC^2$ of the form 
$f(z,w)=(p(z),q(z,w))$ with $f(0,0)=(0,0)$ and $df_{(0,0)}=\mathrm{id}$, 
whose second differential is nondegenerate. Such maps are called 
skew-products tangent to the identity. They showed that $f$ admits a 
wandering domain of rank one under a Diophantine condition on its 
coefficients. 

To give a precise statement, we recall some notation from \cite{Astorg-Boc 
Thaler-2024}. Up to conjugacy by a biholomorphism of $\mathbb{C}^2$, we may 
assume that $f(z,w)$ has the form \cite[(1.1)]{Astorg-Boc Thaler-2024} 
\begin{equation}\label{map1}
\left\{
	\begin{aligned}
	p(z)&:= z-z^2+a z^3+O(|z|^4),\\
    q(z,w)&:=  w+w^2+bz^2+b_{0,3} w^3+b_{3,0} z^3+O(\|(z,w)\|^4),
	\end{aligned}
\right.
\end{equation}
where $a, b, b_{0,3}, b_{3,0} \in \bC$. As in \cite[Theorem 1.6]{Astorg-Boc Thaler-2024}, we assume $b>\frac{1}{4}$ and $b_{0,3}-a\in \R$ and introduce the following parameters: 
\begin{equation}\label{eq:c} 
	c:=\sqrt{b-\tfrac{1}{4}}>0, \quad \alpha:=e^{\pi/c}>1, \quad \beta:=(b_{0,3}-a)(\alpha-1)\in \R.
\end{equation}

\begin{theorem}\label{t:ABT0}\cite[Theorem 1.6]{Astorg-Boc Thaler-2024}
Let $f(z,w)$ be a skew-product of the form~\eqref{map1} with coefficients 
satisfying \eqref{eq:c}. If there exists a strictly increasing sequence of 
positive integers $(n_k)_{k\geqslant 1}$ such that the phase sequence 
$(n_{k+1}-\alpha n_k-\beta\ln n_k)_{k\geqslant 1}$ converges, then $f(z,w)$ 
has wandering domains of rank one. 
\end{theorem}

This theorem builds a bridge between a question in complex dynamics about 
the existence of wandering domains, and a question in number theory about 
the distribution of certain real sequences $(n_{k+1}-\alpha n_k-\beta\ln 
n_k)_{k\geqslant 1}$. 

\begin{question}[Boc Thaler]\label{q:BT}
For which $\alpha>1$ and $\beta\in\mathbb{R}$ does there exist a strictly 
increasing sequence of positive integers $(n_k)_{k\geqslant 1}$ such that 
$(n_{k+1}-\alpha n_k-\beta\ln n_k)_{k\geqslant 1}$ converges? 
\end{question}

Astorg and Boc Thaler proved some almost sharp conditions on 
$(\alpha,\beta)$ for this question. The condition on $\alpha$ is closely 
related to the Pisot property, a notion arising from the study of the 
distribution of powers of $\alpha$ modulo $1$. 

\begin{definition}\label{defn-Pisot-number}\leavevmode
\begin{enumerate}  
\item A real number $\alpha>1$ has the {\em Pisot property} if there exists 
    some non-zero $\xi\in\R$ such that 
    $\lim\limits_{n\to+\infty}\|\xi\alpha^n\|=0$, where $\|\cdot\|$ denotes 
    the distance to the nearest integer. 
\item A \emph{Pisot number} is a real algebraic integer $\alpha>1$ all of 
    whose Galois conjugates other than $\alpha$ have modulus strictly less 
    than $1$. 
\end{enumerate} 
\end{definition}

\begin{theorem}(\cite{Hardy}, \cite{Pisot})
A real algebraic number $\alpha>1$ has the Pisot property if and only if it 
is a Pisot number. 
\end{theorem}

A long-standing open question, dating back to Hardy \cite{Hardy}, asks 
whether there exists any transcendental real number $\alpha>1$ with the Pisot 
property. See \cite[Chapter~2]{Bugeaud-2012} for more on this question. 

Astorg and Boc Thaler considered the following weaker version of convergence. 

\begin{definition}\cite[Definition~1.9]{Astorg-Boc Thaler-2024}
For any positive integer $\ell$,
a sequence $(x_{k})_{k \geqslant 1}$ is said to
{\em converge to a cycle of period $\ell$} 
if the subsequence $(x_{k\ell + j})_{k \geqslant 1}$ converges for each $0 \leqslant j < \ell$.
\end{definition}

Note that the above definition does not require $\ell$ to be the exact 
period. Now we can fully state the conditions of Astorg and Boc Thaler on 
$(\alpha,\beta)$. For the condition on $\beta$ it is convenient to 
introduce the parameter 
$\theta:=\frac{\ln\alpha}{\alpha-1}\beta\in\mathbb{R}$. 

\begin{theorem}\cite[Theorem~1.10]{Astorg-Boc Thaler-2024}\label{t:ABT} Let $\alpha>1$ and let $\beta\in \R$.
\begin{enumerate}
    \item If there exists a strictly increasing sequence of positive 
        integers $(n_k)_{k\geqslant 1}$ such that $(n_{k+1}-\alpha 
        n_k-\beta\ln n_k)_{k\geqslant 1}$ converges to a cycle, then 
        $\alpha$ has the Pisot property. 
    \item Conversely, if $\alpha$ has the Pisot property and  
        $\theta=\frac{\beta\ln\alpha}{\alpha-1}\in \frac{1}{\ell}\Z$, 
        where $\ell\ge 1$ is an integer, then there exists a strictly 
        increasing sequence of positive integers $(n_k)_{k\geqslant 1}$ 
        such that $(n_{k+1}-\alpha n_k-\beta\ln n_k)_{k\geqslant 1}$ 
        converges to a cycle of period $\ell$. 
\end{enumerate}
\end{theorem}

Astorg and Boc Thaler posed the following question on sharp conditions on 
$\beta$ for convergence to a cycle. 

\begin{question}\cite[Question~10.10]{Astorg-Boc Thaler-2024}\label{q:ABT}
Let $\alpha>1$ be a real number with the Pisot property and let $\beta\in 
\R$. Is $\theta:= \frac{\beta\ln\alpha}{\alpha-1}\in \Q$ a necessary 
condition for the existence of a strictly increasing sequence of positive 
integers $(n_k)_{k\geqslant 1}$ such that  $(n_{k+1}-\alpha n_k-\beta\ln 
n_k)_{k\geqslant 1}$ converges to a cycle? 
\end{question}

Assuming that $\alpha$ is algebraic, we answer both Question~\ref{q:BT} and 
Question~\ref{q:ABT} by giving necessary and sufficient conditions on 
$\beta$. In fact, we give necessary and sufficient conditions for each period 
$\ell$. 

\begin{theorem}[Main~Theorem]\label{t:Main}
Let $\alpha>1$ be a real algebraic number with minimal polynomial 
$P(x)\in\mathbb{Z}[x]$. Let $\beta\in \R$ and let $\ell$ be a positive 
integer. There exists a strictly increasing sequence of positive integers 
$(n_k)_{k\geqslant 1}$ such that  $(n_{k+1}-\alpha n_k-\beta\ln 
n_k)_{k\geqslant 1}$ converges to a cycle of period $\ell$ if and only if 
$\alpha$ is a Pisot number and $\theta=\frac{\beta\ln\alpha}{\alpha-1}$ 
belongs to $\frac{1}{\ell P(1)}\mathbb{Z}$. 
\end{theorem}

By the minimal polynomial of $\alpha$, we mean an irreducible polynomial 
$P(x)\in \Z[x]$ such that $P(\alpha)=0$ and the leading coefficient of $P$ 
is positive. In particular, if $\alpha$ is an algebraic integer, then $P$ is 
monic. 

When combined with Theorem \ref{t:ABT0}, the case $\ell=1$ of our Main 
Theorem produces skew-products with wandering domains of rank one.  Note that 
Theorem \ref{t:ABT0} requires the sequence $(n_{k+1}-\alpha n_k-\beta\ln 
n_k)_{k\geqslant 1}$ to converge. It is not yet clear whether $f(z,w)$ still 
has wandering domains when the sequence $(n_{k+1}-\alpha n_k-\beta\ln 
n_k)_{k\geqslant 1}$ converges only to a cycle. After announcing an earlier 
version of this article which answered, for algebraic $\alpha$, 
Question~\ref{q:ABT} on convergence to a cycle, Boc Thaler suggested that we 
consider Question~\ref{q:BT} on convergence.

\begin{corollary}\label{c:wandering}
Let $\alpha>1$ be a real algebraic number with minimal polynomial 
$P(x)\in\mathbb{Z}[x]$. Let $\beta\in \R$. There exists a strictly 
increasing sequence of positive integers $(n_k)_{k\geqslant 1}$ such that  
$(n_{k+1}-\alpha n_k-\beta\ln n_k)_{k\geqslant 1}$ converges if and only if 
$\alpha$ is a Pisot number and $\theta=\frac{\beta\ln\alpha}{\alpha-1}$ 
belongs to $\frac{1}{ P(1)}\mathbb{Z}$. In this case, the skew-product 
$f(z,w)$ of the form \eqref{map1} whose coefficients satisfy \eqref{eq:c} 
has wandering domains of rank one. 
\end{corollary}

\begin{example} Let $\alpha$ be a Pisot number with minimal polynomial $P(x)\in \Z[x]$. 
For any $j\in\mathbb{Z}$, the polynomial skew-product $f(z,w)=(p(z),q(z,w))$ 
defined by 
\begin{equation}\label{eq:ex}
\left\{\begin{aligned}
p(z)&:=z-z^2,\\
q(z,w)&:=w+w^2+\frac{j}{P(1)\ln\alpha}w^3+\left(\frac{1}{4}+\frac{\pi^2}{(\ln \alpha)^2}\right) z^2,
\end{aligned}\right.
\end{equation}
admits wandering domains of rank one by Corollary \ref{c:wandering} applied 
to $\theta=j/P(1)$. The case $P(1)\nmid j$ is new. In particular, every 
integer $\alpha>2$ is a Pisot number with $P(1)=1-\alpha$ for which 
\eqref{eq:ex} gives explicit new examples by choosing $j$ not divisible by 
$\alpha-1$. 
\end{example}

\begin{example} 
Let $\alpha>1$ be an integer. Then the minimal polynomial $P(x)=x-\alpha$ 
satisfies $P(1)=1-\alpha$ and the condition in Corollary \ref{c:wandering} 
becomes $\beta\ln \alpha \in \Z$. Thus, for any $\beta'\in \R\backslash 
\Z$, there exist no strictly increasing sequences of positive integers 
$(n_k)_{k\geqslant 1}$ such that the sequence $(n_{k+1}-\alpha 
n_k-\beta'\log_\alpha n_k)_{k\geqslant 1}$ converges. 
\end{example}

For the sufficiency in the Main Theorem, we construct the sequences 
explicitly using complete homogeneous polynomials. Using linear recurrence 
relations and Weyl's equidistribution theorem, we prove that the necessity of 
the rationality of $\theta$ in the Main Theorem holds more generally, for 
example under the assumption that the phase sequence is bounded and the set 
of limit values is countable.

This paper is organized as follows. In Section~\ref{s:2}, we prove the Main 
Theorem using complete homogeneous polynomials. In Section~\ref{s:3}, we 
prove a generalization of the necessity for rationality in the Main Theorem. 

\section{Proof of the Main Theorem}\label{s:2}

Before proving the Main Theorem, we will prove the following related result 
on the limit of $(m_{k+1}-\alpha m_k)_{k\ge 1}$. 

\begin{prop}\label{p:limit}
Let $\alpha>1$ be a Pisot number and let $P(x)\in \Z[x]$ be its minimal 
polynomial. Let $\gamma\in \R$. Then the following conditions are 
equivalent: 
\begin{enumerate}[(a)]
\item There exists a strictly increasing sequence of integers 
    $(m_k)_{k\ge 1}$ such that $(m_{k+1}-\alpha m_k)_{k\ge 1}$ converges 
    to $\gamma$. 
\item There exists a sequence of integers $(m_k)_{k\ge 1}$ such that 
    $(m_{k+1}-\alpha m_k)_{k\ge 1}$ converges to $\gamma$. 
\item $\gamma\in \frac{1-\alpha}{P(1)}\cdot\Z$. 
\end{enumerate}
\end{prop}

The key to proving the existence of the sequence of integers in both the 
Main Theorem and Proposition \ref{p:limit} is the use of complete 
homogeneous polynomials. We let 
\[h_k(x_1,\dots,x_d)=\sum_{1\le i_1\le \dots\le i_k\le d} x_{i_1}\dotsm x_{i_k}\] 
denote the complete homogeneous polynomial of degree 
$k$. 

\begin{lemma}\label{lem}
Let $\alpha_1,\dots,\alpha_d,\zeta\in \CC$ such that $\lvert \alpha_j\rvert 
<1$ for all $1\leqslant j\leqslant d$ and $\lvert \zeta\rvert = 1$. Let 
$Q(x)=(1-\alpha_1 x)\dotsm (1-\alpha_d x)$. Then, as $k\to \infty$, we have 
\begin{enumerate}
\item[(a)] $h_k(1,\alpha_1,\dots,\alpha_d)=\frac{1}{Q(1)}+o(1)$.
\item[(b)] 
    $h_k(1,\zeta,\alpha_1,\dots,\alpha_d)
    =\begin{cases}\frac{k+1}{Q(1)}+\frac{Q'(1)}{Q(1)^2}+o(1)&\zeta=1\\
    \frac{1}{1-\zeta}\cdot\left(\frac{1}{Q(1)}-\frac{\zeta^{k+1}}{Q(\zeta^{-1})}\right)+o(1)&\zeta\neq 
    1.\end{cases}$
\end{enumerate}
\end{lemma}

\begin{proof}
(a) Note that
\[
h_k(1,\alpha_1,\dots,\alpha_d)=\sum\limits_{t=0}^k 
h_t(\alpha_1,\dots,\alpha_d)
\]
is the sum of the coefficients of the terms of degrees $\leqslant k$ of the 
power series
\begin{equation}\label{eq:Qinv} 
\frac{1}{Q(x)}=\sum\limits_{t\geqslant 0} h_t(\alpha_1,\dots,\alpha_d)x^t.
\end{equation}
Since $|\alpha_j|< 1$ for each $1\leqslant j\leqslant d$, the power series 
converges absolutely when $|x|\leqslant 1$. Thus $\lim\limits_{k\to \infty} 
h_k(1,\alpha_1,\dots,\alpha_d)=\frac{1}{Q(1)}$. 

(b) For $\zeta=1$, 
\[
h_k(1,1,\alpha_1,\dots,\alpha_d)=
\sum\limits_{t=0}^k h_t(1,\alpha_1,\dots,\alpha_d)
=
\sum\limits_{t=0}^{k}(k+1-t)h_t(\alpha_1,\dots,\alpha_d)
\]
is the sum of the coefficients of the terms of degrees $\leqslant k$ of the power 
series
\[
\frac{k+1}{Q(x)}+x\frac{Q'(x)}{Q(x)^2}
=
\sum\limits_{t\geqslant 0}(k+1-t)h_t(\alpha_1,\dots,\alpha_d)x^t.
\]
Since $|\alpha_j|< 1$ for all $1\leqslant j\leqslant d$, the sum of the 
coefficients of the terms of degrees $>k$ is $o(1)$ as $k\to \infty$.  Thus 
$h_k(1,1,\alpha_1,\dots,\alpha_d)=\frac{k+1}{Q(1)}+\frac{Q'(1)}{Q(1)^2}+o(1)$. 
For $\zeta\neq 1$, 
\begin{align*}
h_k(1,\zeta,\alpha_1,\dots,\alpha_d)&=\sum_{t=0}^k \sum_{s=0}^{k-t}\zeta^{s}h_t(\alpha_1,\dots,\alpha_d)=\sum_{t=0}^k\frac{1-\zeta^{k-t+1}}{1-\zeta}\cdot h_t(\alpha_1,\dots,\alpha_d)\\
&=\frac{1}{1-\zeta}\cdot\left(\frac{1}{Q(1)}-\frac{\zeta^{k+1}}{Q(\zeta^{-1})}\right)+o(1)
\end{align*}
by \eqref{eq:Qinv}.
\end{proof}

\begin{lemma}\label{lem2}
Let $\alpha_1,\dots,\alpha_d\in \CC$ be stable under complex conjugation. 
Assume that $\alpha_1$ is real and $\alpha_1 
>\lvert \alpha_j\rvert$  for all $1< j\leqslant d$. Then 
$h_k(\alpha_1,\dots,\alpha_d)/\alpha_1^k$ converges to a positive real number 
as $k\to \infty$. 
\end{lemma}

\begin{proof}
By homogeneity, 
\[
h_k(\alpha_1,\dots,\alpha_d)/\alpha_1^k=h_k(1,\tfrac{\alpha_2}{\alpha_1},\dots,\tfrac{\alpha_d}{\alpha_1}),
\]
which converges to $1/\prod\limits_{j=2}^d (1-\alpha_j/\alpha_1)>0$ by Lemma 
\ref{lem}(a). 
\end{proof}

\begin{proof}[Proof of Proposition \ref{p:limit}]
(a)$\implies$(b). Trivial.

(b)$\implies$(c). Let $d$ be the degree of $P(x)$ and let $Q(x)=x^dP(1/x)$ 
and $R(x)=Q(x)/(1-\alpha x)=\sum_{s=0}^{d-1} b_s x^s$. Consider the formal 
power series $M=\sum_{k=1}^\infty m_k x^k$. By assumption, the coefficient 
of $x^k$ in $(1-\alpha x)M$ converges to $\gamma$ as $k\to \infty$. Thus 
the coefficient of $x^k$ in $R(x)(1-\alpha x)M=Q(x)M\in \Z[[x]]$ converges 
to $(b_0+\dots+b_{d-1})\gamma\in \Z$. Since 
$b_0+\dots+b_{d-1}=R(1)=\frac{P(1)}{1-\alpha}$, it follows that 
$\frac{P(1)}{1-\alpha}\cdot \gamma\in \Z$. 

(c)$\implies$(a). Let $\alpha=\alpha_1,\alpha_2,\dots,\alpha_d$ be the 
Galois conjugates of $\alpha$. Let $w_k=h_k(1,\alpha_1,\dots,\alpha_d)$. 
Since $w_k$ is an algebraic integer in $\Q$, we have $w_k\in \Z$. Moreover,
\[w_{k+1}-\alpha w_k=h_{k+1}(1,\alpha_2,\dots,\alpha_d)\to \frac{1}{(1-\alpha_2)\dotsm (1-\alpha_d)}=\frac{1-\alpha}{P(1)}\]
as $k\to \infty$ by Lemma \ref{lem}(a) and the assumption that $\alpha$ is 
a Pisot number. We have $\gamma=(1-\alpha)(\frac{a}{P(1)}+b)$ for some 
integers $a,b$ with $a>0$. Take $m_k=aw_k+b$. By Lemma \ref{lem2}, $m_k$ is 
strictly increasing for sufficiently large $k$. It then suffices to discard 
finitely many initial terms of $(m_k)_{k\ge 1}$. 
\end{proof}

The proof of the necessity in the Main Theorem uses the following result by 
Astorg and Boc Thaler. 

\begin{lemma}\label{l:cycle}\cite[Lemma~10.4]{Astorg-Boc Thaler-2024}
Let $\alpha>1$ and $\beta\in \R$. Let $(n_k)_{k\geqslant 1}$ be a strictly 
increasing sequence of positive integers, and let 
$\ell\in\mathbb{Z}_{\geqslant 1}$. Define $m_k:=n_{k+\ell}-n_k$. Assume 
that the sequence $(n_{k+1}-\alpha n_k-\beta\ln n_k)_{k\geqslant 1}$ 
converges to a cycle of period $\ell$. Then the sequence $(m_{k+1}-\alpha 
m_k)_{k\geqslant 1}$ converges to $\ell \beta \ln \alpha$. 
\end{lemma}

\begin{proof}[Proof of Theorem \ref{t:Main}]
(Necessity) Suppose there is a strictly increasing sequence of positive 
integers $(n_k)_{k\geqslant 1}$ such that $(n_{k+1}-\alpha n_k-\beta\ln 
n_k)_{k\geqslant 1}$ converges to a cycle of period $\ell$. Let 
$m_k:=n_{k+\ell}-n_k$. Then $(m_{k+1}-\alpha m_k)_{k\geqslant 1}$ converges 
to $\ell \beta\ln\alpha$ by Lemma~\ref{l:cycle}. By Theorem \ref{t:ABT}(1), 
$\alpha$ has the Pisot property and hence is a Pisot number. By Proposition 
\ref{p:limit}, it follows that $\ell \beta\ln\alpha\in 
\frac{1-\alpha}{P(1)}\Z$. 

(Sufficiency) Let $\alpha_1=\alpha,\alpha_2,\dots,\alpha_d$ be the Galois 
conjugates of $\alpha$. Let 
\[
w_k:=\sum_{j=0}^{\lfloor\frac{k}{\ell}\rfloor} h_{k-j\ell}(1,\alpha_1,...,\alpha_d)\in\mathbb{Z}.
\]
Since $\frac{\beta\ln\alpha}{\alpha-1}\in\frac{1}{\ell P(1)}\mathbb{Z}$, 
there exist some integers $a,b$ with $a>0$ such that 
$\frac{\beta\ln\alpha}{1-\alpha}=\frac{a}{\ell P(1)}+b$. Define 
$n_k:=aw_k+bk$. For any integer $s$, 
\[\sum_\zeta \zeta^s =\begin{cases}
\ell & \ell\mid s\\
0 & \ell\nmid s,
\end{cases}\] 
where $\zeta$ runs through all $\ell$-th roots of unity. Thus
\[\ell w_k=\sum_{\zeta}\sum_{s=0}^k \zeta^s h_{k-s}(1,\alpha_1,\dots,\alpha_d)=\sum_{\zeta} h_{k}(1,\zeta,\alpha_1,\dots,\alpha_d).\]
It follows that
\begin{equation}\label{eq:wdiff}
w_{k+1}-\alpha w_{k}=\frac{1}{\ell}\sum_\zeta h_{k+1}(1,\zeta,\alpha_2,\dots,\alpha_d)=\frac{k}{\ell(1-\alpha_2)\dotsm(1-\alpha_d)}+C_{\bar k}+o(1)
\end{equation}
for some constant $C_{\bar k}$ depending on the congruence class $\bar k$ 
of $k$ modulo $\ell$ by Lemma \ref{lem}(b) and the fact that $\alpha$ is a 
Pisot number. Moreover, since 
\[\prod_\zeta \frac{1}{1-\zeta x}=\frac{1}{1-x^\ell},\]
we have
\[w_k=h_k(1,\alpha_1,\dots,\alpha_d,\zeta_1,\dots,\zeta_\ell),\]
where $\zeta_1,\dots,\zeta_\ell$ denote the $\ell$-th roots of unity. Thus, 
by Lemma~\ref{lem2}, $n_k>0$ for $k$ sufficiently large and 
\begin{equation}\label{eq:lnnk}
\ln n_k = k\ln \alpha + C'+o(1)
\end{equation}
for some constant $C'$. In particular, $n_k$ is strictly increasing for $k$ 
sufficiently large. By \eqref{eq:wdiff} and \eqref{eq:lnnk}, 
\[
\aligned
n_{k+1}-\alpha n_k -\beta\ln n_k
&=\frac{a(1-\alpha)k}{\ell P(1)}+a C_{\bar k}+(1-\alpha) b k+b-\beta k\ln \alpha -\beta C'+o(1)\\
&=a C_{\bar k}+b-\beta C'+o(1)
\endaligned
\]
converges to a cycle of period $\ell$.
\end{proof}

\section{A generalization}\label{s:3}

In this section, we prove a generalization of the necessity of the 
rationality of $\theta$ in the Main Theorem, replacing convergence to a cycle 
by certain restrictions on the distribution of the phase sequence. To give a 
precise statement, we introduce the following notion. 

\begin{defn}\label{d:thin}
Let $(a_0,\dots,a_d)$ be a finite collection of real numbers. We say that a 
bounded subset $I\subseteq \R$ is \emph{$(a_0,\dots,a_d)$-thin} if the 
closure of the Minkowski sum $\sum_{s=0}^d a_s I=\{\sum_{s=0}^d a_s x_s\mid 
x_0,\dots,x_d\in I\}\subseteq \R$ has Lebesgue measure zero. We say that a 
bounded subset $I\subseteq \R$ is \emph{thin} if it is $(a_0,\dots,a_d)$-thin 
for every  finite collection of real numbers $(a_0,\dots,a_d)$.
\end{defn}

\begin{example}\label{ex:thin}
Any compact subset $I\subseteq \R$ of lower Minkowski dimension $<1/(d+1)$ 
or of packing dimension $<1/(d+1)$ is $(a_0,\dots,a_d)$-thin. In 
particular, any countable compact subset $I\subseteq\R$ is thin. We refer 
the reader to \cite{Falconer} for the definition and properties of these 
dimensions. Packing dimension was originally introduced by Tricot 
\cite{Tricot}. 
\end{example}

Given a sequence $(x_k)_{k\geqslant 1}$ and a set $J$, denote by 
\[
A(x_\bullet,J):=\{k\geqslant 1\mid x_k\in J\}.
\]

\begin{definition}
For any set $A\subseteq \mathbb{Z}_{\geqslant 1}$, define its {\em upper 
uniform density} and {\em lower uniform density} as 
\begin{align*}
\uden(A)&:=\lim_{n\to \infty}\frac{1}{n}\sup_{h\in \Z_{\ge 0}}\#(A\cap [h+1,h+n]),\\
\lden(A)&:=\lim_{n\to \infty}\frac{1}{n}\inf_{h\in \Z_{\ge 0}}\#(A\cap [h+1,h+n]).
\end{align*}
We say that $A$ has \emph{uniform density} $\delta$ if  
$\uden(A)=\lden(A)=\delta$. 
\end{definition}

We refer to \cite{GTT} for a proof that the limits exist and for equivalent 
definitions. 

\begin{theorem}\label{t:nec}
Let $\alpha>1$ be a real algebraic number and let $P(x)=a_d x^d+\dotsm 
+a_0\in \R[x]$ be a polynomial such that $P(1)\neq 0$ and $P(\alpha')=0$ 
for every Galois conjugate $\alpha'\neq \alpha$ of~$\alpha$. Let $\beta\in 
\R$. Let $(n_k)_{k\ge 1}$ be a strictly increasing sequence of positive 
integers. Assume that $(\sigma_k)_ {k\geqslant 1}:=(n_{k+1}-\alpha 
n_k-\beta \ln n_k)_{k\geqslant 1}$ is bounded and there exists an 
$(a_0,\dots,a_d)$-thin subset $I\subseteq \R$ such that $\delta:=\inf_J 
\uden(A(\sigma_\bullet,J))>1-\frac{1}{d+1}$, where $J$ runs through open 
neighborhoods of~$I$. Then $\theta:=\frac{\beta\ln\alpha}{\alpha-1}\in \Q$. 
\end{theorem}

\begin{remark}
The condition in the theorem is satisfied if $(\sigma_k)_{k\geqslant 1}$ is 
bounded and the set $E$ of limit values of the sequence is 
$(a_0,\dots,a_d)$-thin. Indeed, for any open neighborhood $J$ of $E$, $\{k\ge 
1\mid \sigma_k \not\in J\}$ is finite. 

If a sequence converges to a cycle, then the set of limit values is finite 
and, in particular, thin. 
\end{remark}

Our proof of Theorem \ref{t:nec} is different from that of the necessity 
part of the Main Theorem given in Section \ref{s:2}. We recall some 
constructions and a result from Astorg and Boc Thaler \cite[Section 
10]{Astorg-Boc Thaler-2024} which will be used in the proof. Let 
$\alpha>1$, $\beta\in\mathbb{R}$ and 
$\theta:=\frac{\beta\ln\alpha}{\alpha-1}$. For any strictly increasing 
sequence of positive integers $(n_k)_{k\geqslant 1}$, define 
\begin{align*}
\sigma_k&:=n_{k+1}-\alpha n_k-\beta\ln n_k,\\
m_k&:=n_k+\lfloor k\theta\rfloor,\\
\eta_k&:=m_{k+1}-\alpha m_k+\{(k+1)\theta\}-\alpha\{k\theta\}.
\end{align*}

\begin{lemma}\label{l:ABT}
Assume that $(\sigma_k)_{k\ge 1}$ is bounded. Then 
$\sigma_k=\eta_k+\lambda+o(1)$ for some constant $\lambda\in \R$. In 
particular, for any closed subset $I\subseteq \R$ and any open neighborhood 
$J$ of $I$, there are only finitely many $k\ge 1$ such that $\sigma_k\in I$ 
but $\eta_k+\lambda\not\in J$. 
\end{lemma}

\begin{proof} 
The first assertion is proved in \cite[Remark~10.9]{Astorg-Boc 
Thaler-2024}. For completeness, let us recall the proof. By the second 
displayed formula in \cite[Lemma~10.1]{Astorg-Boc Thaler-2024}, there 
exists a real number $\xi>0$ such that 
\begin{align}\label{e:sigma-eta}
\sigma_k&=n_{k+1}-\alpha n_k-\beta(k\ln\alpha+\ln\xi)+o(1)\\
\notag&=m_{k+1}-\alpha m_k+\{(k+1)\theta\}-\alpha\{k\theta\}-\theta-\beta\ln\xi+o(1)\\
\notag&=\eta_k-\theta-\beta\ln\xi+o(1).
\end{align}
The last assertion follows from the first one. 
\end{proof}

By Lemma \ref{l:ABT}, it suffices to study the distribution of 
$(\eta_k)_{k\ge 1}$. Our key idea here is to study the distribution of the 
linear combination
\begin{equation}\label{d:tildek}
\tilde{\eta}_k:=\sum\limits_{s=0}^d a_s \eta_{k+s}.
\end{equation}
By definition, we have
\begin{equation}\label{ekck}
\tilde{\eta}_k=e_{k} + c_k,
\end{equation}
where
\begin{align}\label{d:ek}
e_k&:=\sum\limits_{s=0}^d a_{s} (m_{k+s+1}-\alpha m_{k+s}),\\
\label{d:ck}
c_k&:=a_d\{(k+d+1)\theta\}
+
\left(\sum\limits_{s=1}^{d}(a_{s-1}-\alpha a_s)\{(k+s)\theta\}\right)
-\alpha a_0 \{k\theta\}.
\end{align}

\begin{lemma}\label{L:ek}
Assume that $(\sigma_k)_{k\geqslant 1}$ is bounded. Then $(e_k)_{k\geqslant 
1}$ defined by \eqref{d:ek} takes only finitely many values. 
\end{lemma}

\begin{proof} 
Let $Q(x)=\sum_{i=0}^g u_ix^i\in \Z[x]$ be the minimal polynomial of 
$\alpha$. Then $(x-\alpha)P(x)=Q(x)R(x)$ for some $R(x)=\sum_{i=0}^t 
v_ix^i\in \R[x]$. Let $f_k=\sum_{i=0}^g u_i m_{k+i}\in \Z$. By \cite[Lemma 
10.1]{Astorg-Boc Thaler-2024}, there exists $\xi>0$ such that 
$(m_k-\xi\alpha^k)_{k\ge 1}$ is bounded. Thus 
\[f_k=\sum_{i=0}^g u_i (m_{k+i}-\xi\alpha^{k+i})+\xi\alpha^k\underbrace{\sum_{i=0}^g u_i \alpha^i}_{=Q(\alpha)=0}\]
is bounded. Consequently $(f_k)_{k\ge 1}$ takes only finitely many values. It 
follows that the same holds for $(e_k)_{k\ge 1}$, since $e_k=\sum_{i=0}^t 
v_if_{k+i}$. 
\end{proof} 

When $\theta\notin \mathbb{Q}$, the distribution of the sequence 
$(c_k)_{k\geqslant 1}$ can be characterized using irrational rotations. 
Define a piecewise linear function 
\[
\aligned
g:[0,1) &\longrightarrow \mathbb{R}\\
x&
\longmapsto
a_d\{x+(d+1)\theta\}
+
\left(\sum\limits_{s=1}^{d}(a_{s-1}-\alpha a_s)\{x+s\theta\}\right)
-\alpha a_0 x,
\endaligned
\]
so that $c_k=g(\{k\theta\})$. The set of discontinuity points of $g$ is a 
subset of 
\[\bigl\{ \{-\theta\}, \{-2\theta\}, \dots, \{-(d+1)\theta\} 
\bigr\},\] 
where the elements are irrational and mutually distinct. These 
$d+1$ points partition the interval $(0, 1)$ into $d+2$ open intervals 
$J_0, J_1, \dots, J_{d+1}$. Define $I_j := g(J_j)$ for $j = 0, 1, \dots, 
d+1$. 

Let $\Omega$ be a compact topological space and let $\mu$ be a regular 
positive Borel measure on $\Omega$ such that $\mu(\Omega)=1$. A sequence 
$(x_k)_{k\ge 1}$ in $\Omega$ is said to be \emph{$\mu$-well distributed} if 
for every continuous function $\varphi$ on $\Omega$, we have $\lim_{n\to 
\infty}\frac{1}{n}\sum_{k=h+1}^{h+n}\varphi(x_k)= \langle 
\varphi,\mu\rangle$ uniformly in $h\ge 0$. We refer to \cite[Chapter~3, 
Corollaries 3.1, 3.2]{KN} for equivalent definitions. In particular, by 
\cite[Chapter~3, Corollary 3.2]{KN}, $(x_k)_{k\ge 1}$ is $\mu$-well 
distributed if and only if for every Borel set $M\subseteq \Omega$ 
satisfying $\mu(\partial M)=0$, $A(x_\bullet,M)$ has uniform density equal 
to $\mu(M)$. The notion of well-distribution was introduced by Hlawka 
\cite{Hlawka0} and Petersen \cite{Petersen}. 

\begin{lemma}\label{l:Leb}
Assume $\theta \notin \mathbb{Q}$. Let $\Omega\subseteq \R$ be a compact 
interval containing $g([0,1))$ and let 
\[\mu=\frac{1}{(\alpha-1)\lvert P(1)\rvert}\sum\limits_{j=0}^{d+1}\chi_{I_{j}} \mleb_\Omega,\]
where $\mleb_\Omega$ denotes the Lebesgue measure on $\Omega$ and 
$\chi_{I_j}$ is the characteristic function of $I_j$. Then $\mu(\Omega)=1$ 
and $(c_k)_{k\ge 1}$ is $\mu$-well distributed.
\end{lemma}

\begin{proof} 
Since $\theta\notin \bQ$, by an extension \cite{Petersen} (see also 
\cite[Chapter~1, Example 5.2]{KN}) of Weyl's equidistribution theorem 
\cite[Satz~2]{Weyl}, $(\{k\theta\})_{k\ge 1}$ is $\mleb_{[0,1]}$-well 
distributed. The lemma follows easily by pushing forward via $g$. In more 
detail, for any continuous function $\varphi$ on $\Omega$, we have 
\begin{equation}\label{eq:leb} 
\lim_{n \to \infty} \frac{1}{n} \sum_{k=h+1}^{h+n} \varphi ( c_k )  = 
\lim_{n \to \infty}\frac{1}{n} \sum_{k=h+1}^{h+n} \varphi\left(g ( \{ k\theta   \}   ) \right)
=\int_{0}^{1}\varphi\left(g(x)\right)dx = \sum\limits_{j=0}^{d+1} \int_{J_j} 
\varphi\left(g(x)\right)dx 
\end{equation}
uniformly in $h\ge 0$, by \cite[Theorem~2]{Petersen} applied to the 
Riemann-integrable function $\varphi\circ g$. For each $j=0,1,\dots,d+1$, 
$g$ is a linear function on $J_j$ with derivative 
\begin{align*}
g'(x)&=a_d+\left(\sum\limits_{s=1}^{d}(a_{s-1}-\alpha a_s)\right)-\alpha a_0\\
&=(a_d+a_0+\dots+a_{d-1})-\alpha(a_0+a_1+\dots+a_d)=(1-\alpha)P(1)
\end{align*}
which is a constant independent of $j$. Since $P(1)\neq 0$ by assumption, 
we conclude that $g'(x)\neq 0$ for any $x\in J_j$. Thus 
\begin{equation}\label{eq:leb2}
\sum\limits_{j=0}^{d+1} \int_{J_j} \varphi\left(g(x)\right)dx = 
\sum\limits_{j=0}^{d+1} \int_{I_j}\frac{1}{(\alpha-1)\lvert P(1)\rvert} \varphi(y)dy=
\langle\varphi,\mu\rangle.
\end{equation}
Combining \eqref{eq:leb} and \eqref{eq:leb2}, we get
\[\lim_{n \to \infty} \frac{1}{n} \sum_{k=h+1}^{h+n} \varphi ( c_k )  = \langle\varphi,\mu\rangle\]
uniformly in $h\ge 0$. In particular, taking $\varphi=1$, we get 
$\mu(\Omega)=1$. 
\end{proof}

Our next lemma is about lower uniform density. In general lower uniform 
density fails to be subadditive: There exist subsets $A,B\subseteq \Z_{\ge 
1}$ such that $\lden(A\cup B)>\lden(A)+ \lden(B)$. However, we have the 
following lemma for the lower uniform density of a finite union of 
translates. For $A\subseteq \Z_{\ge 1}$ and $m\in \Z$, we put $A[m]=\{k\ge 
1\mid k+m\in A\}$. 

\begin{lemma}\label{l:upper-den}
Let $A\subseteq\mathbb{Z}_{\geqslant1}$ and let $M$ be a finite set of 
integers. Then $\lden(\bigcup_{m\in M}A[m]) \le \#M\cdot \lden(A)$ and 
$1-\uden(\bigcap_{m\in M}A[m]) \le \#M\cdot (1- \uden(A))$. 
\end{lemma}

\begin{proof}
We may assume $M$ nonempty. Let $a_n=\#( A\cap [1,n])$. For $n\le 0$ we 
adopt the convention that $a_n=0$. Then $a_n\le a_{n+1}\le a_n+1$ for all 
$n\in \Z$. Let $U=\bigcup_{m\in M}A[m]$ and $k=\max_{m\in M} \lvert 
m\rvert$. Then 
\begin{multline*}
\#(U\cap 
[h+1,h+n])\le \sum_{m\in M}\#(A[m]\cap[h+1,h+n])=\sum_{m\in M}(a_{h+n+m}-a_{h+m})\\\le \#M\cdot (a_{h+n+k}-a_{h-k})\le \#M\cdot 
(a_{h+n}-a_h+2k)=\#M\cdot (a_{h+n}-a_h) +O(1).
\end{multline*}
Thus 
\[\lden(U)=\lim_{n\to \infty} \frac{1}{n}\inf_{h\ge 0}\#(U\cap 
[h+1,h+n])\le\#M\cdot \lim_{n\to \infty}\frac{1}{n}\inf_{h\ge 0} (a_{h+n}-a_h)=\#M\cdot \lden(A).
\] 
The inequality for upper uniform density follows by taking complement. 
Indeed, $\Z_{\ge 0}\backslash (A[m])$ and $(\Z_{\ge 0}\backslash A)[m]$ 
differ only by finitely many elements. 
\end{proof} 

\begin{proof}[Proof of Theorem~\ref{t:nec}] 
Suppose $\theta\notin\mathbb{Q}$. Up to replacing $I$ by its closure, we 
may assume that $I$ is closed. By assumption, $S=\sum_{s=0}^d a_s I$ has 
Lebesgue measure zero. Let $\epsilon>0$. There exists a union $U\subseteq 
\R$ of finitely many open intervals such that $S\subseteq U$ and 
$\mleb(U)<\epsilon$, where $\mleb$ denotes the Lebesgue measure. By 
continuity, there exists an open neighborhood $J$ of $I$ such that 
$\sum_{s=0}^d a_s J\subseteq U$. Let $J_0\subseteq J$ be a closed 
neighborhood of $I$. By assumption, $\uden(A(\sigma_\bullet, J_0)) \ge 
\delta$. Thus, by Lemma \ref{l:ABT}, $\uden(A(\eta_\bullet+\lambda, J)) \ge 
\delta$. By \eqref{ekck}, 
\[
c_k=\tilde{\eta}_{k}-e_{k}.
\]
Let $E=\{-e_{k}\mid k\ge 1\}$, which is a finite set by Lemma \ref{L:ek}. 
Then 
\begin{equation}\label{eq:Aincl}
A\big(c_\bullet+P(1)\lambda,U+E\big)\supseteq A\big(\tilde{\eta}_\bullet+P(1)\lambda,U\big)\supseteq \bigcap_{s=0}^d A(\eta_{\bullet+s}+\lambda,J).
\end{equation}
By Lemma~\ref{l:upper-den}, the last set in \eqref{eq:Aincl} has upper 
uniform density $\ge1-(d+1)(1-\delta)>0$, and thus so does the first set 
$A\big(c_\bullet+P(1)\lambda,U+E\big)$. Let $\mu$ be the measure in Lemma 
\ref{l:Leb}, which satisfies $\mu\le C\mleb_\Omega$ for some constant 
$C>0$. By \cite[Chapter 3, Corollary 3.2]{KN}, since $(U+E-P(1)\lambda)\cap 
\Omega$ is a $\mu$-continuity set, $A\big(c_\bullet+P(1)\lambda,U+E\big)$ 
has uniform density equal to $\mu\big((U+E-P(1)\lambda)\cap \Omega\big)\le 
C\mleb\big(U+E-P(1)\lambda\big)< C \epsilon \cdot \# E$. Taking 
$\epsilon>0$ sufficiently small produces the desired contradiction. 
\end{proof} 

\subsection*{Acknowledgements} Zhangchi Chen was supported by National Key 
R\&D Program of China (No.2025YFA1018300), NSFC (grant number 12501104), 
the Science and Technology Commission of Shanghai Municipality (Grant No. 
22DZ2229014), the Shanghai Sailing Program (Grant No. 24YF2709900), the 
Shanghai Pujiang Program (Grant No. 24PJA023). Zhangchi Chen also 
acknowledges support from the Labex CEMPI (ANR-11-LABX-0007-01), the 
project QuaSiDy (ANR-21-CE40-0016) and the CDP C2EMPI 
(R-CDP-24-004-$C^2$EMPI). Weizhe Zheng was supported by the National 
Natural Science Foundation of China (grant numbers 12125107, 12288201) and 
the Chinese Academy of Sciences Project for Young Scientists in Basic 
Research (grant number YSBR-033). 

We sincerely thank Matthieu Astorg and Luka Boc Thaler for proposing this 
question and for explaining their key ideas from their work. We extend our 
special thanks to Luka Boc Thaler for comments on an early version of our 
manuscript and for suggesting that we consider Question~\ref{q:BT} in a 
private conversation. We are grateful to Zhuchao Ji for hosting the first 
and second authors at Westlake University and for discussions. We thank 
Charles Favre for hosting the first author at {\'E}cole Polytechnique and 
suggesting that we read the paper by Astorg and Boc Thaler. We thank 
Shaoshi Chen, Wenxia Li, Xiao Yao, and Longtu Yuan for inspiring 
discussions. We thank Wenjuan Peng and Fei Yang for sharing their knowledge 
on dynamical systems. We thank Weikun He for the suggestion to consider the 
packing dimension in Example \ref{ex:thin}.


\begin{bibdiv}
\begin{biblist}
\bib{ABDPR-2016}{article}{
   author={Astorg, Matthieu},
   author={Buff, Xavier},
   author={Dujardin, Romain},
   author={Peters, Han},
   author={Raissy, Jasmin},
   title={A two-dimensional polynomial mapping with a wandering Fatou
   component},
   journal={Ann. of Math. (2)},
   volume={184},
   date={2016},
   number={1},
   pages={263--313},
   issn={0003-486X},
   review={\MR{3505180}},
   doi={10.4007/annals.2016.184.1.2},
}

\bib{Astorg-Boc Thaler-2024}{article}{
   author={Astorg, Matthieu},
   author={Boc Thaler, Luka},
   title={Dynamics of skew-products tangent to the identity},
   journal={J. Eur. Math. Soc. (JEMS)},
   volume={28},
   date={2026},
   number={2},
   pages={559--618},
   issn={1435-9855},
   review={\MR{5031425}},
   doi={10.4171/jems/1566},
}

\bib{Bugeaud-2012}{book}{
   author={Bugeaud, Yann},
   title={Distribution modulo one and Diophantine approximation},
   series={Cambridge Tracts in Mathematics},
   volume={193},
   publisher={Cambridge University Press, Cambridge},
   date={2012},
   pages={xvi+300},
   isbn={978-0-521-11169-0},
   review={\MR{2953186}},
   doi={10.1017/CBO9781139017732},
}

\bib{Falconer}{book}{
   author={Falconer, Kenneth},
   title={Fractal geometry},
   edition={3},
   note={Mathematical foundations and applications},
   publisher={John Wiley \& Sons, Ltd., Chichester},
   date={2014},
   pages={xxx+368},
   isbn={978-1-119-94239-9},
   review={\MR{3236784}},
}

\bib{GTT}{article}{ author={Grekos, Georges}, author={Toma, Vladim\'ir}, 
   author={Tomanov\'a, Jana}, title={A note on uniform or Banach density}, 
   journal={Ann. Math. Blaise Pascal}, volume={17}, date={2010}, 
   number={1}, pages={153--163}, issn={1259-1734}, review={\MR{2674656}}, 
} 

\bib{Hardy}{article}{
   author={Hardy, Godfrey Harold},
   title={A problem of Diophantine approximation},
   journal={J. Indian Math. Soc.},
   volume={11},
   date={1919},
   pages={162--166},
}

\bib{Hlawka0}{article}{
   author={Hlawka, Edmund},
   title={Zur formalen Theorie der Gleichverteilung in kompakten Gruppen},
   language={German},
   journal={Rend. Circ. Mat. Palermo (2)},
   volume={4},
   date={1955},
   pages={33--47},
   issn={0009-725X},
   review={\MR{0074489}},
   doi={10.1007/BF02846027},
}

\bib{KN}{book}{
   author={Kuipers, L.},
   author={Niederreiter, H.},
   title={Uniform distribution of sequences},
   series={Pure and Applied Mathematics},
   publisher={Wiley-Interscience [John Wiley \& Sons], New
   York-London-Sydney},
   date={1974},
   pages={xiv+390},
   review={\MR{0419394}},
}

\bib{Petersen}{article}{
   author={Petersen, G. M.},
   title={`Almost convergence'\ and uniformly distributed sequences},
   journal={Quart. J. Math. Oxford Ser. (2)},
   volume={7},
   date={1956},
   pages={188--191},
   issn={0033-5606},
   review={\MR{0095812}},
   doi={10.1093/qmath/7.1.188},
}

\bib{Pisot}{article}{
   author={Pisot, Charles},
   title={La r\'epartition modulo 1 et les nombres alg\'ebriques},
   language={French},
   journal={Ann. Scuola Norm. Super. Pisa Cl. Sci. (2)},
   volume={7},
   date={1938},
   number={3-4},
   pages={205--248},
   issn={0391-173X},
   review={\MR{1556807}},
}

\bib{Sullivan-1985}{article}{
   author={Sullivan, Dennis},
   title={Quasiconformal homeomorphisms and dynamics. I. Solution of the
   Fatou-Julia problem on wandering domains},
   journal={Ann. of Math. (2)},
   volume={122},
   date={1985},
   number={3},
   pages={401--418},
   issn={0003-486X},
   review={\MR{0819553}},
   doi={10.2307/1971308},
}

\bib{Tricot}{article}{ author={Tricot, Claude, Jr.}, title={Two definitions 
   of fractional dimension}, journal={Math. Proc. Cambridge Philos. Soc.}, 
   volume={91}, date={1982}, number={1}, pages={57--74}, issn={0305-0041}, 
   review={\MR{0633256}}, doi={10.1017/S0305004100059119}, 
} 

\bib{Weyl}{article}{
   author={Weyl, Hermann},
   title={\"Uber die Gleichverteilung von Zahlen mod. Eins},
   language={German},
   journal={Math. Ann.},
   volume={77},
   date={1916},
   number={3},
   pages={313--352},
   issn={0025-5831},
   review={\MR{1511862}},
   doi={10.1007/BF01475864},
}
\end{biblist}
\end{bibdiv}
\end{document}